\documentclass[12pt]{amsart}
\usepackage[T1]{fontenc}
\usepackage{lmodern}
\usepackage{amsmath,amssymb}
\usepackage[a4paper,textwidth=13.5cm,textheight=24cm]{geometry}
\usepackage{microtype}
\usepackage{xurl}
\usepackage[hidelinks]{hyperref}
\hypersetup{
  pdftitle={On integers whose friable component is large},
  pdfauthor={Gerald Tenenbaum and Xiangyu Wang},
  pdfsubject={Uniform estimates for integers with a large friable component},
  pdfkeywords={friable integers, smooth numbers, Dickman function, saddle-point method}
}
\newcommand{\e}{{\rm e}}
\newtheorem{theorem}{Theorem}

\title[Integers whose friable component is large]{On integers whose friable\\component is large}
\author[G.~Tenenbaum]{G\'erald Tenenbaum}
\address{Institut \'Elie Cartan, Universit\'e de Lorraine, BP 70239,
54506 Vand\oe uvre Cedex, France}
\email{gerald.tenenbaum@univ-lorraine.fr}
\author[X.~Wang]{Xiangyu Wang}
\address{Department of Mathematics, University of Illinois Urbana-Champaign,
Urbana, Illinois, USA}
\email{xw70@illinois.edu}
\subjclass[2020]{11N25}
\keywords{Friable integers, smooth numbers, Dickman function, saddle-point method}
\date{}

\begin{document}

\begin{abstract}
We estimate the number of integers up to $x$ whose $y$-friable component
exceeds $z$. Uniformly for $\max\{2,\log z\}\leqslant  y\leqslant  z\leqslant  x/2$, we
determine the logarithm of this counting function up to an error
$O(\log z/\log y)$. This sharpens an estimate arising from a problem in
Ford's sieve methods lecture notes. The proof combines classical
friable-number estimates with an elementary argument for the lower bound.
\end{abstract}

\maketitle

\section{Introduction}

For $y\geqslant 2$, the $y$-friable component $n_y$ of a positive integer $n$
is its largest divisor all of whose prime factors are at most $y$. Put
\[
\Theta(x,y,z):=\#\{n\leqslant  x:n_y>z\},
\qquad v:=\frac{\log z}{\log y}.
\]
Estimates for this function involving the Dickman and Buchstab functions
were obtained in~\cite{Tenenbaum2006}.

Our starting point is Exercise~5.2 of Ford's sieve methods
lecture notes~\cite[p.~59]{Ford2023}, which asks for the estimate
\[
\Theta(x,y,z)=xv^{-v}\e^{O(v\log_2(3v))}
\]
uniformly for $10\leqslant \log z\leqslant  y\leqslant  z\leqslant  x/2$. Here and throughout,
$\log_2 t=\log\log t$. The second author obtained a proof while studying
this problem and subsequently presented it in~\cite{Wang2026}.
We sharpen this estimate by a short application of known results.

\begin{theorem}\label{thm:main}
Uniformly for $\max\{2,\log z\}\leqslant  y\leqslant  z\leqslant  x/2$, we have
\[
\Theta(x,y,z)
=\frac{x\e^{O(v)}}{\{v\log(3v)\}^{v}}.
\]
The implied constant is absolute.
\end{theorem}

\paragraph{Notation.}
We write $P^+(n)$ and $P^-(n)$ for the largest and smallest prime factors
of $n$, with the convention $P^+(1)=1$, $P^-(1)=\infty$, and put
\[
\Psi(t,y):=\#\{n\leqslant  t:P^+(n)\leqslant  y\},
\qquad
\Phi(t,y):=\#\{n\leqslant  t:P^-(n)>y\}.
\]
The letters $\varrho$ and $\omega$ denote the Dickman and Buchstab functions,
respectively. All sums over integers are over positive integers unless
otherwise specified; implied constants are absolute.

\section{Proof of Theorem~\ref{thm:main}}\label{sec:proof}

Put $L=\log y$ and $F(v)=v\log\{v\log(3v)\}$. We may assume that $y$
is sufficiently large. Indeed, bounded $y$ implies bounded $z$, since
$z\leqslant  \e^y$. Choose a power of $2$, say $q$, with $z<q\leqslant 2z$.
Every multiple of $q$ is counted, so
$\Theta(x,y,z)\geqslant \lfloor x/q\rfloor\geqslant  x/4z\gg x$.
Together with $\Theta(x,y,z)\leqslant  x$, this proves the assertion in that case.

For $x>  zy$, decomposition into friable and sifted parts gives
\[
\Theta(x,y,z)=
\sum_{\substack{z<m\leqslant  x/y\\P^+(m)\leqslant  y}}
\Phi(x/m,y)+\Psi(x,y)-\Psi(x/y,y).
\]
For $x\leqslant zy$, we instead have exactly
$\Theta(x,y,z)=\Psi(x,y)-\Psi(z,y)$.
We use the standard bounds $\Phi(t,y)\ll t/L$ for $t\geqslant  y$ and
$\Phi(t,y)\gg t/L$ for $t\geqslant 2y$.

\subsection{\texorpdfstring{The Dickman range: $y>\e^{(\log_2 z)^2}$}{The Dickman range}}

In this range $\log(vL)<\sqrt L$. Classical estimates
\cite[Theorem~1.1 and Corollary~2.3]{HildebrandTenenbaum1993},
\cite[(1.2), (1.8)]{Tenenbaum2006} give
\[
\varrho(v)=\e^{-F(v)+O(v)},
\qquad
\tau(v):=\int_v^\infty\varrho(t)\,dt=\e^{-F(v)+O(v)},
\]
\[
S(y,z):=\sum_{\substack{m>z\\P^+(m)\leqslant  y}}\frac1m
\ll \tau(v)L,
\qquad
\Psi(t,y)=t\varrho(v_t)\{1+o(1)\},
\]
where $z\leqslant  t\leqslant  zy^2$ and $v_t=(\log t)/L$.
The latter $o(1)$ is uniform as $y\to\infty$ in the present range.
Estimates for ratios of Dickman functions yield
$\Psi(2z,y)/\Psi(z,y)=2+o(1)$ and
$\varrho(v+h)=\varrho(v)\e^{O(\log(v+2))}$ for $0\leqslant  h\leqslant 2$.
Consequently, if $2z\leqslant  x\leqslant  zy^2$, then
\[
\Theta(x,y,z)\geqslant \Psi(x,y)-\Psi(z,y)
\gg\Psi(x,y)=x\e^{-F(v)+O(v)}.
\]
For the corresponding upper bound, use $\Theta(x,y,z)\leqslant \Psi(x,y)$ when
$x\leqslant  zy$; when $zy<x\leqslant  zy^2$, the preceding decomposition gives
\[
\Theta(x,y,z)\ll\frac{xS(y,z)}L+\Psi(x,y)
\ll x\e^{-F(v)+O(v)}.
\]

If $x>zy^2$, put $u:=(\log x)/L$. Formula~(1.7) of the corrected version
of~\cite{Tenenbaum2006} applies and yields
\[
\Theta(x,y,z)\asymp x\vartheta(u,v),
\qquad
\vartheta(u,v):=\varrho(u)+\int_v^{u-1}\varrho(t)\omega(u-t)\,dt.
\]
Indeed, its hypotheses hold with, for example, $\varepsilon=1/20$, since
$\log_2 z<L^{1/2}$, $z\geqslant  y$, and $u>v+2$.
Since $1/2\leqslant \omega(s)\leqslant 1$ for $s\geqslant 1$ and
$\int_v^{v+1}\varrho(t)\,dt\asymp\tau(v)$, we obtain
$\vartheta(u,v)\asymp\tau(v)$. This proves both bounds in the Dickman range.

\subsection{\texorpdfstring{The complementary range: $y\leqslant \e^{(\log_2 z)^2}$}{The complementary range}}

Here $v\geqslant  \e^{\sqrt L}/L$, so $L=o(v)$, and the hypothesis $\log z\leqslant  y$
gives $1\leqslant  v\leqslant  y/L$. Let $\beta=\alpha(z,y)$ be the saddle point,
defined by
\[
\sum_{p\leqslant  y}\frac{\log p}{p^\beta-1}=\log z=vL,
\]
and write
\[
\zeta(s,y):=\prod_{p\leqslant  y}(1-p^{-s})^{-1},
\qquad Z:=z^{\beta-1}\zeta(\beta,y),
\qquad w:=(1-\beta)L.
\]
Theorem~1.6 and Lemmas~4.1--4.2 of~\cite{HildebrandTenenbaum1993}
(see also~\cite{HildebrandTenenbaum1986}) give
\[
\frac{c}{L}\leqslant \beta<1,
\qquad \frac{\e^w-1}{w}\asymp v,
\qquad w=\log\{v\log(3v)\}+O(1),
\]
for an absolute $c>0$, together with
\[
\frac{\Psi(z,y)}z\asymp\frac{Z}{\beta L\sqrt v}.
\]
To estimate the Euler product, Chebyshev's bounds imply, uniformly for
$\beta\leqslant \sigma\leqslant 1$,
\[
\sum_{p\leqslant  y}\frac{\log p}{p^\sigma-1}
\ll L\left\{1+\frac{\e^{(1-\sigma)L}-1}{(1-\sigma)L}\right\},
\]
with the quotient interpreted by continuity at $\sigma=1$.
Integrating and using $\zeta(1,y)\asymp L$, we find
\[
0\leqslant \log\zeta(\beta,y)
\ll\log_2(3y)+w+\int_0^w\frac{\e^t-1}{t}\,dt\ll v.
\]
Thus $\log Z=-vw+O(v)=-F(v)+O(v)$. The saddle-point formula and
$L=o(v)$ now give
\[
Z=\e^{-F(v)+O(v)},
\qquad \frac{\Psi(z,y)}z=Z\e^{O(v)}.
\]
Rankin's bound gives $S(y,z)\leqslant  Z$ and
$\Psi(x,y)\leqslant  x^\beta\zeta(\beta,y)\leqslant  xZ$.
Using the decomposition above, or $\Theta\leqslant \Psi(x,y)$ when $x\leqslant y$,
we conclude that $\Theta(x,y,z)\ll xZ$.

For the lower bound, map each $y$-friable integer $m\leqslant  z$ to $2^jm$,
where $j$ is the least non-negative integer such that $2^jm>z$.
The image belongs to $(z,2z]$ and remains $y$-friable. Each image has at
most $1+\log(2z)/\log2$ preimages. Hence
\[
D:=\Psi(2z,y)-\Psi(z,y)\gg\frac{\Psi(z,y)}{\log(2z)}.
\]
If $x>4zy$, restrict the friable part to $(z,2z]$ and use
$\Phi(x/m,y)\gg x/mL$. If $2z\leqslant  x\leqslant 4zy$, simply use
$\Theta(x,y,z)\geqslant  D$. These two arguments give, respectively,
\[
\Theta(x,y,z)\gg
\begin{cases}
\displaystyle\frac{x\Psi(z,y)}{zL\log(2z)},&x>4zy,\\[6pt]
\displaystyle\frac{x\Psi(z,y)}{yz\log(2z)},&2z\leqslant  x\leqslant 4zy.
\end{cases}
\]
Since $L=o(v)$ and $\log_2(2z)\ll L$, the extra factors in both denominators
have logarithm $O(v)$. The preceding estimates for $Z$ and $\Psi(z,y)$
therefore yield $\Theta(x,y,z)\geqslant  x\e^{-F(v)-Cv}$ for an absolute $C>0$.
This completes the proof.\hfill$\square$

\subsection*{Acknowledgments}
The second author thanks Kevin Ford for reading an earlier version of his
manuscript.

\subsection*{Use of generative AI}
The second author used ChatGPT to assist with English drafting, LaTeX preparation, reference checking, and the exposition of existing mathematical arguments.

\end{document}